\documentclass[3p,preprint,nopreprintline]{elsarticle}

\usepackage{lineno}

\usepackage{hyperref}
\hypersetup{pdfauthor={Name}}

\usepackage{amsfonts}
\usepackage{graphicx}
\usepackage{commath, mathtools}
\usepackage{subfig}
\usepackage{booktabs}
\usepackage{soul}
\usepackage{xcolor}
\usepackage{placeins}

\usepackage{multirow}

\usepackage[capitalise,noabbrev]{cleveref}
\crefformat{equation}{(#2#1#3)}

\usepackage[algo2e,ruled,vlined,linesnumbered]{algorithm2e}
\SetCommentSty{scriptsize}
\SetKwComment{tcc}{[}{]}
\SetKwIF{If}{ElseIf}{Else}{if}{}{else if}{else}{endif}
\SetKwFor{While}{while}{}{endw}
\SetKwFor{ForEach}{for each}{}{endfch}
\SetKwFor{For}{for}{}{}%
\SetKwComment{tcp}{$\rightarrow$\quad}{}
\SetNlSty{}{\color{gray}}{\quad}
\SetKwInOut{Input}{input}
\SetKwInOut{Output}{output}

\Crefname{algocf}{Algorithm}{Algorithms}

\def\ddefloop#1{\ifx\ddefloop#1\else\ddef{#1}\expandafter\ddefloop\fi}
\def\ddef#1{\expandafter\def\csname bf#1\endcsname{\ensuremath{\mathbf{#1}}}}
\ddefloop ABCDEFGHIJKLMNOPQRSTUVWXYZabcdefghijklmnopqrstuvwxyz\ddefloop

\def\ddef#1{\expandafter\def\csname v#1\endcsname{\ensuremath{\boldsymbol{#1}}}}
\ddefloop ABCDEFGHIJKLMNOPQRSTUVWXYZabcdefghijklmnopqrstuvwxyz\ddefloop

\def\ddef#1{\expandafter\def\csname v#1\endcsname{\ensuremath{\boldsymbol{\csname #1\endcsname}}}}
\ddefloop {alpha}{beta}{gamma}{delta}{epsilon}{varepsilon}{zeta}{eta}{theta}{vartheta}{iota}{kappa}{lambda}{mu}{nu}{xi}{pi}{varpi}{rho}{varrho}{sigma}{varsigma}{tau}{upsilon}{phi}{varphi}{chi}{psi}{omega}{Gamma}{Delta}{Theta}{Lambda}{Xi}{Pi}{Sigma}{varSigma}{Upsilon}{Phi}{Psi}{Omega}{ell}\ddefloop

\def\ddef#1{\expandafter\def\csname bb#1\endcsname{\ensuremath{\mathbb{#1}}}}
\ddefloop ABCDEFGHIJKLMNOPQRSTUVWXYZ\ddefloop

\usepackage{amsopn}

\DeclareMathOperator{\bspan}{span}
\DeclareMathOperator{\argmin}{arg\,min}
\newcommand{\rb}{\text{rb}}

\newcommand{\train}{\text{train}}

\newcommand{\Dtrain}{\mathcal{D}_{\train}}
\newcommand{\Ntrain}{N_{\train}}

\begin{document}

\begin{frontmatter}

\title{Reduced Basis Approximations of Parameterized Dynamical Partial Differential Equations via Neural Networks\tnoteref{t1}}

\author[2]{Peter Sentz\corref{cor1}}
\ead{sentz2@illinois.edu}

\author[1]{Kristian Beckwith}
\ead{kbeckwi@sandia.gov}

\author[1]{Eric C. Cyr}
\ead{eccyr@sandia.gov}

\author[2]{Luke N. Olson}
\ead{lukeo@illinois.edu}
\ead[url]{http://lukeo.cs.illinois.edu}

\author[1]{Ravi Patel}
\ead{rgpatel@sandia.gov}

\affiliation[1]{%
  organization={Sandia National Laboratories},
  city={Albuquerque},
  state={NM}}
\affiliation[2]{%
  organization={University of Illinois at Urbana-Champaign},
  addressline={Department of Computer Science},
  city={Urbana},
  state={IL}}
\cortext[cor1]{Corresponding author}

\begin{abstract}
 Projection-based reduced order models are effective at approximating parameter-dependent differential equations that are parametrically separable.  When parametric separability is not satisfied, which occurs in both linear and nonlinear problems, projection-based methods fail to adequately reduce the computational complexity.  Devising alternative reduced order models is crucial for obtaining efficient and accurate approximations to expensive high-fidelity models.  In this work, we develop a time-stepping procedure for dynamical parameter-dependent problems, in which a neural-network is trained to propagate the coefficients of a reduced basis expansion.  This results in an online stage with a computational cost independent of the size of the underlying problem.  We demonstrate our method on several parabolic partial differential equations, including a problem that is not parametrically separable.
\end{abstract}


\end{frontmatter}

\section{Introduction}

In this paper, we develop a method for obtaining reduced order models (ROMs) of time-dependent parameterized partial differential equations (PDEs). Our method uses proper orthogonal decomposition (POD) to obtain reduced space representations of the system. It learns dynamics on this reduced space, thereby avoiding the need to project the full-order model in the online stage and resulting in a runtime that is independent of the full-order model size.

Complex physical phenomena and their underlying models often include a wide
variety and number of parameters, from material properties to geometric factors
to initial conditions.  In parameter studies, such as uncertainty
quantification or design optimization, it is often computationally infeasible to
simulate each parameter combination at full resolution.  Instead, the modeling can be
decomposed into \textit{offline} and \textit{online} stages in order to reduce
cost.  In the offline stage, high resolution approximations or so-called full-order
solutions are found for a small sample of $P$ parameters, from which a basis is constructed. The online stage uses projection to
construct solutions for new parameter combination.

The use of projection-based reduced order modeling is well understood for
stationary PDEs~\cite{rozza2007reduced,chaudyry2021lsrbm,PACCIARINI20161977}, while
the treatment for time-based schemes is less developed, particularly when the
dynamics of the problem evolve rapidly.  One recent
development~\cite{HESTHAVEN201855} is the use of feed-forward neural networks
to construct projection coefficients; in essence, the problem parameters are
mapped to the coefficients.  Extensions of this approach to time-based
problems~\cite{WaHeRa_2019_nonintrusivenn} include the time value as an additional input feature, again yielding a map
directly to the coefficients.

The approach we propose in this work is to
construct a mapping from coefficients at time step $k$ to time step $k+1$, i.e., pairs of coefficients at different time steps are used to develop or \textit{learn}
a mapping from one state to the next. Thus the feed-forward network is required only to predict the mapping from one time step to the next. The previously mentioned neural network approach, in contrast, must learn and encode the entire dynamics over the time range of interest. While this may lend benefit to the proposed approach, it is important to consider that a discrete time stepping error is also incurred.
Thus we view our technique as an alternative to the prior neural network methodology. Further, compared to using a Galerkin projection in the online phase, the computational cost of the proposed method does not scale with the size of the full order problem in the context of non-affine parameter dependence. This computational scaling improvement and the temporal iteration aspect of our technique, have the potential to expand the range of dynamics accessible to ROM approaches.

Other works consider error models~\cite{Freno_2019, parish2020time} to handle the dynamics
in time.  The error is modeled as a random variable, and recurrent neural
networks learn a prediction model in time.  Methods also extend to the
\textit{training} portion~\cite{fresca2020comprehensive}, where the focus is on
learning an effective trial space and the mapping over time is still
monolithic, as in~\cite{WaHeRa_2019_nonintrusivenn}.  In addition, other
methods~\cite{xu2019multi}, use a multilayered approach for handling time and
space (through the use of autoencoders), while still remaining non-intrusive
(where the PDE itself is not used in the construction).

The paper and contributions are organized as follows.  In ~\cref{sec:param_pdes}, we outline the basics of parameterized PDEs and
the mechanics of reduced-order models using reduced basis methods and
projected-based schemes.  In~\cref{sec:romres} we introduce an approach for
computing projection coefficients using neural network regression.  We motivate
the use of a residual neural network (ResNet) and develop the full algorithm
in~\cref{sec:algorithm}.  In~\cref{sec:numerics} we consider several time
dependent problems; we highlight the efficacy of the method developed in this
work and discuss several salient choices in the algorithm.

\section{Parametrized PDEs and Reduced Order Modeling}\label{sec:param_pdes}

In this work, we consider time-dependent partial differential equations (PDEs) of the form:
\begin{equation}\label{eq:general_pde}
		u_t+ F(u; \vmu) = 0,\quad(\vx,t) \in \Omega\times (0,T],
\end{equation}
where $u = u(\vx, t; \vmu)$ is the unknown solution, $F(u;\vmu)$ is a differential operator with respect to the spatial variables, and $\Omega\subset\mathbb{R}^d$ is the spatial simulation domain.  In addition, $F$ depends on a vector of $P$ parameters $\vmu = (\mu_1,\dots,\mu_P)^T \in \mathcal{D}\subset \mathbb{R}^P$, that represent physical or geometric constants.  For example, a linear advection equation with wave speed $\mu$ leads to the operator
\begin{equation}
	F(u; \mu) = \mu u_x,
\end{equation}
with parameter dimension $P = 1$.  In other situations, the differential operator $F$ is nonlinear with multiple parameters.  For example, the parametrized porous medium equation~\cite{drohmann2012reduced} features the operator
\begin{equation}
		F(u; \vmu) = -\mu_1\Delta(u^{\mu_2}),
\end{equation}
where $\vmu = (\mu_1, \mu_2) \in \mathbb{R}^+\times\mathbb{R}$.  \ref{sec:notation} summarizes the notation used throughout the remainder of the manuscript.

For many problems, the solution algorithm demands that a parametrized PDE be solved repeatedly (the ``many-query'' context)%
~\cite{rozza2007reduced}.  For example, in a PDE-constrained optimization problem, a PDE is solved at every step of the optimization loop.  In other cases, a real-time measurement of a physical quantity demands a rapid solution of a PDE for prediction of future behavior.  A high-fidelity (or \textit{full-order}) numerical approximation to~\cref{eq:general_pde} is often prohibitively expensive in either of these contexts.  This leads to so-called \textit{reduced order} modeling: introducing methods that replace a numerical scheme with one of much lower complexity and thus reduced computation time, while retaining an acceptable level of accuracy.

\subsection{Reduced Basis Methods}

A commonly used family of reduced-order models are \textit{reduced basis methods}~\cite{quarteroni2015reduced}, which are based on a Galerkin projection.  Generically, a Galerkin projection approximates the solution to~\cref{eq:general_pde} by a function of the form
\begin{equation}\label{eq:fe_solution}
		u^h(\vx,t;\vmu) \coloneqq \sum_{j=1}^{N_h} \alpha_j(t;\vmu)\phi_j(\vx),
\end{equation}
where $\{\phi_1(\vx),\dots,\phi_{N_h}(\vx)\}$ is a linearly independent set of functions in the spatial variable, with time/parameter-dependent coefficients $\valpha(t;\vmu) = [\alpha_1(t;\vmu),\dots,\alpha_{N_h}(t;\vmu)]^T$.  Here, $h > 0$ is a parameter that determines the size of the basis; e.g., the size of a grid in a finite element or finite difference discretization.

The coefficients in~\cref{eq:fe_solution} are solved for using a variational formulation;~\cref{eq:general_pde} is multiplied by $\phi_i \in V^h\coloneqq \text{span}\{\phi_1,\dots,\phi_{N_h}\}$ and integrated over the domain $\Omega$.
This leads to the projected system of ordinary differential equations (ODE):
\begin{equation}\label{eq:general_ode}
		\sum_{j=1}^{N_h}\alpha_j'(t;\vmu)\int_{\Omega}\phi_j\phi_i\ \dif \vx + \int_{\Omega}F\left(\sum_{j=1}^{N_h} \alpha_j(t;\vmu)\phi_j(\vx);\vmu\right)\phi_i\ \dif \vx = 0,
\end{equation}
for $1 \leq i \leq N_h$.
This ODE is solved by an appropriate time-stepping method.  In the context of reduced order modeling,~\cref{eq:general_ode} is known as the \emph{full-order} model.

Instead of projecting onto $V^h$, a reduced basis method projects the problem in~\cref{eq:general_pde} onto a space of $N_{\rb}$ \textit{reduced basis} functions $V^{\rb} =
\bspan\{\psi_1(\vx),\dots, \psi_{N_{\rb}}(\vx)\}$, yielding a semi-discrete reduced model.  A temporal discretization can then be applied either before or after the introduction of the reduced basis.

For computational efficiency, we seek a reduced basis with $\dim(V^{\rb}) = N_{\rb} \ll N_h$.  However, this by itself does not guarantee a satisfactory time complexity that is independent of the full-order dimension $N_h$.  The key to a successful reduced basis approach is the separation of the algorithm into a (possibly expensive) \textit{offline} stage, and a faster \textit{online} stage. In the \textit{offline} stage the reduced basis is constructed to resolve a set of high-fidelity solutions.
The approach we use for constructing the basis is the Proper Orthogonal Decomposition (POD), that is briefly reviewed in the next section.
The runtime cost of the \textit{offline} stage typically scales with the size of the full-order problem.
In contrast, the \textit{online} stage ideally has computational complexity that is independent of the original problem size.
In this stage,  the projection of the reduced-order solution of~\cref{eq:general_ode} is computed for novel parameter instances. %

\subsection{Offline Stage: Construction of the Basis}\label{sec:pod_construction}

There are numerous ways to select a basis of smaller dimension to represent solutions to PDEs.
In this paper, we use the POD algorithm to obtain $V^{\rb}$~---~the neural network model learns the time evolution of the POD coefficients for a corresponding parameter $\vmu$.  Consequently, the approach is applicable to any method in which the reduced-order model is represented as a linear combination of basis functions with $N_{\rb} \ll N_h$.
Moreover, while POD is often used with a variety of approximation schemes~---~finite volume, finite difference, etc.~---~we restrict the presentation here to finite elements, although our approach does not rely on this.

In the standard POD approach the offline stage builds a matrix $A$ of solution snapshots at sampled parameter values and time steps. As shown in~\cref{alg:pod}, the POD basis is constructed by computing the SVD of $A$, and truncating to a specified tolerance denoted by $\epsilon$.
Each column of the resulting matrix, $\hat{A}$, corresponds to the degrees of freedom of a reduced basis function $\psi_k(x)$. We employ one additional step: through the use of the mass matrix, we impose full $L^2$ orthogonality of the basis, rather than discrete~$\ell^2$ orthogonality, yielding an optimal projection onto the discrete function space.  Additional details can be found in~\cite{quarteroni2015reduced}.

The traditional process is effective for constructing a basis, but can be expensive if the number of parameters and time steps is large.
To alleviate this, we follow~\cite{WaHeRa_2019_nonintrusivenn} and adopt a two-stage process to construct the basis, which results in a theoretically sub-optimal, but computationally efficient alternative to the standard POD algorithm.
In the first stage, a finite parameter sample $\mathcal{D}_{\text{POD}} = \{\vmu_1, \vmu_2,\dots,\vmu_{N_\text{s}}\} \subset \mathcal{D}$ of $N_\text{s}$ samples is constructed.  For each parameter sample $\vmu_j$, with $j\in \{1,\dots, N_\text{s}\}$, a matrix containing the degrees of freedom of the full-order solutions is constructed:
\begin{equation}
	S_j = \left[\valpha(t_1;\vmu_j)\ \ \valpha(t_2;\vmu_j)\ \  \dots\ \ \valpha(t_{N_t};\vmu_j) \right].
\end{equation}
An SVD, ${S}_j = \hat{U}_j\Sigma_j V_j^T$, is formed, and the first $N_{j}$ columns of $\hat{U}_j$ are retained, denoted by $U_j$.
This corresponds to a compression over time for a fixed parameter $\vmu$.  A standard selection criteria for the number of columns $N_j$ is the smallest integer $m$ such that
\begin{equation}\label{eq:pod-criterion}
	\frac{\sum_{i = m+1}^{r}\sigma_i^2}{\sum_{i=1}^{r}\sigma_i^2}\leq \epsilon_{t},
\end{equation}
where $r$ is the rank of $S_j$, $\sigma_i$ are singular values, and $\epsilon_t$ is a parameter controlling the amount of compression in time.

In the second stage, the matrices ${U}_1,\dots,{U}_{N_\text{s}}$ from the first stage are concatenated into a matrix
\begin{equation}
  {U} = \left[{U}_1,\dots,{U}_{N_{\text{s}}}\right].
\end{equation}
An SVD is again performed to obtain ${U} = \hat{{W}}{D}{Z}^T$, and the first $N_{\rb}$ columns of $\hat{{W}}$ are retained, to obtain the matrix ${W} \in \mathbb{R}^{N_h\times N_{\rb}}$, using selection criteria~\cref{eq:pod-criterion} with parameter $\epsilon_{\vmu}$.  This second stage corresponds to a compression over the parameter sample $\mathcal{D}_{\text{POD}}$.

\subsection{Online Stage: Computing the POD coefficients}

The online stage assumes the offline stage has computed a basis of much smaller dimension, $N_{\rb} \ll N_h$.  This by itself will not guarantee a sufficient reduction in computational complexity.  The online stage should in fact be \textit{independent} of the problem size $N_h$~\cite{rozza2007reduced}.

If the underlying PDE is linear, then the ODE in~\cref{eq:general_ode} is also linear, and applying an implicit single-step time discretization leads to a sequence of linear algebraic equations at time step $k$ of the form:
\begin{equation}\label{eq:single_step_implicit}
	A(\vmu)\valpha(t_{k+1};\vmu) = B(\vmu)\valpha(t_k;\vmu).
\end{equation}
Using the Galerkin method to project~\cref{eq:single_step_implicit} onto the reduced basis leads to a sequence of problems of much smaller dimension:
\begin{align}
  {W}^T A(\vmu){W}\vc(t_{k+1};\vmu) & = {W}^T B(\vmu){W}\vc(t_k;\vmu),\\
  \implies \vc(t_{k+1};\vmu) &= \left({W}^T A(\vmu){W}\right)^{-1}{W}^T B(\vmu){W}\vc(t_k;\vmu) \label{eq:pod-update}
\end{align}
where $\vc(t_k;\vmu) \in \mathbb{R}^{N_{\rb}}$ is the vector of reduced basis coefficients at time $t_k$.  While this has replaced the need to solve $N_h\times N_h$ systems of linear equations by those of dimension $N_{\rb}\times N_{\rb}$, the assembly of the matrices ${W}^T A(\vmu) {W}$ and ${W}^T B(\vmu) {W}$ for new parameter values still incurs a cost dependent on the problem size $N_h$.

Efficient $N_h$-independent online stages can be implemented in the case that the system has the property of \textit{parametric separability}, i.e.
\begin{align}\label{eq:parameter_separable}
  \begin{split}
		A(\vmu) &= \sum_{j}\theta_j^A(\vmu)A_j,\\
		B(\vmu) &= \sum_{j}\theta_j^B(\vmu)B_j,
  \end{split}
\end{align}
where $A_j$, $B_j$ are parameter independent $N_{\rb}\times N_{\rb}$ matrices, and $\theta_j^A$, $\theta_j^B$ are smooth, real-valued functions of the parameter $\vmu$.  In this case, the assembly of the reduced matrices can be performed with a computational cost independent of problem size in the online stage.  For example, consider ${W}^T A(\vmu){W}$, and use the parametric separability~\cref{eq:parameter_separable} to obtain
\begin{equation}
	{W}^T A(\vmu) {W} = \sum_{j}\theta_j^A(\vmu){W}^T A_j{W}.
\end{equation}
The matrices ${W}^T A_j{W}$ can be computed in the offline stage (at a $N_h$-dependent cost) and stored, while in the online stage, only the scalar-valued functions $\theta_j^A(\vmu)$ must be computed for new parameter values; the assembly of ${W}^T A(\vmu) {W}$ is completed by forming a linear combination of $N_{\rb}\times N_{\rb}$ matrices.  A similar construction is performed for ${W}^T B(\vmu) {W}$.

While there are extensions to problems that do not have parametric separability~\cite{chaturantabut2010nonlinear,drohmann2012reduced}, the end-result is often highly non-trivial, especially for nonlinear problems, and may require substantial changes to computational code~\cite{HESTHAVEN201855}.  Thus, alternatives to projection must be explored to obtain efficient online stages for arbitrary problems.

\section{Time-stepping ROM with Residual Neural Networks}\label{sec:romres}

The ROM approach introduced in the previous section has been successfully applied in many different scenarios.  Yet, a remaining challenge is to develop algorithms for the online phase for general nonlinear problems that scales independently with size of the full-order systems.  To address this, we develop an approach that is motivated in part by the work in~\cite{WaHeRa_2019_nonintrusivenn},  where a neural network is trained to represent a mapping $\mathcal{M}: (t_k,\vmu) \rightarrow \vc(t_k;\vmu)$.
Instead, we approach the problem in a fundamentally different way: by learning a map between the coefficients at consecutive time-steps, namely $\mathcal{N}: (\vc(t_k;\vmu),\vmu) \rightarrow \vc(t_{k+1};\vmu)$.  While the parameter vector $\vmu$ is still an argument of the map, the first argument is now the value of the coefficients at time $t_k$, rather than the value of time itself.  In this section we propose a method to construct the mapping $\mathcal{N}$.

\subsection{Neural Network Regression to Compute POD Coefficients}

The POD algorithm of~\cref{sec:param_pdes} results in a reduced order representation of the form
\begin{equation}
	u^{\text{rb}}(\vx, t;\vmu) \coloneqq \sum_{j=1}^{N_{\text{rb}}}c_j(t;\vmu)\psi_j(\vx).
\end{equation}
Here, the functions $\psi_j(\vx)$ are the reduced basis~---~i.e.\ the columns of matrix $W$.
The coefficients $c_j(t;\vmu)$ are the reduced coefficients, which are the components of a vector $\vc(t;\vmu) \in \mathbb{R}^{N_{\text{rb}}}$.

The overarching goal for the online phase is to construct a mapping from the triple $(u(\vx,0;\vmu),t,\vmu)$ to coefficients $\vc(t;\vmu)$.
A direct approach to constructing this mapping, as described in~\cref{sec:param_pdes}, uses a time-discrete representation of $u^{\text{rb}}$ and
a Galerkin projection of the full-order model.  The update equation~\cref{eq:pod-update} can be used repeatedly, following a projection of the initial condition onto the reduced basis.  However, the computationally efficient application of this method is limited to linear PDEs that are parametrically separable.

An alternative approach is regression with a nonlinearly parameterized model for mappings involving the coefficients.
For instance, in~\cite{WaHeRa_2019_nonintrusivenn} a parameterized model of the form (indices are excluded for brevity)
\begin{equation}\label{eq:timemap}
	\mathcal{N}_R:(t,\vmu) \to \vc(t;\vmu)
\end{equation}
is used.  To construct the model, a training data set is defined by a collection of features (parameters and time points)
\begin{equation}\label{eq:training-data}
\mathcal{F} = \{(k,\vmu_i): k=1\ldots N_t, i=1\ldots N_s \},
\end{equation}
along with corresponding coefficient targets
\begin{equation}
	\mathcal{T} = \{\vc(t_k;\vmu_i) : k=1\ldots N_t, i=1\ldots N_s \}.
\end{equation}
For simplicity, we have described the features and targets using the same set of parameters $\mathcal{D}_{\text{POD}} = \{\vmu_1,\dots,\vmu_{N_s}\}$ and time steps $k=1,\dots,N_t$ that were used to construct the POD basis.  This is not a necessary restriction, and a separate parameter/time sampling from the POD construction may be used.
Using this data set, a loss function is defined on a minibatch $\mathcal{B} \subset\mathcal{F}$ as
\begin{equation}\label{eq:timemap-loss}
		L_{R}(\mathcal{N}_R;\mathcal{B}) = \frac{1}{|\mathcal{B}|} \sum_{(k,\vmu)\in \mathcal{B}} \left\|\mathcal{N}_R(t_k,\vmu) - \vc(t_k,\vmu)\right\|_2^2.
\end{equation}
Here the inputs to the map $\mathcal{N}_R$, $\{(k,\vmu)\}$, and the outputs, $\{\vc(t,\vmu_k)\}$, are determined from
the projection of snapshots used to construct the reduced basis.  In this framework, the optimal parameterizaton of the mapping is constructed by minimizing the loss defined on the full data set $\mathcal{F}$, e.g. $\mathcal{N}_R = \argmin_{\mathcal{N}} L_R(\mathcal{N};\mathcal{F})$.

\subsection{Time-stepping regression method}\label{sec:time-stepping-regression}

Constructing the neural network to define the map in~\cref{eq:timemap} is an ambitious task.  For rapidly evolving dynamics, the difference in POD coefficients for times $t_1$ and $t_1 + \epsilon$ can be relatively large.  Consequently, training data must include a large number of time samples in regions of rapid dynamics. %

An alternative approach that we introduce in this paper is to train a neural network to take time steps, rather than solve a regression problem.  Specifically, we seek to learn the mapping
\begin{equation}\label{eq:update-timemap}
	\mathcal{N}_T:(\vc(t_k;\vmu),\vmu) \to \vc(t_{k+1};\vmu).
\end{equation}
Compared to \cref{eq:timemap}, the primary difference is in the inclusion of the previous coefficients in the domain of the mapping.  Additionally, the mapping is required to learn a \textit{transition} between two POD coefficient states, rather than state prediction at an arbitrary time in the future.  We take advantage of this structure in~\cref{sec:rnns} by learning only the update to the previous coefficient.  This architecture reflects the nature of memoryless ODEs that dependent only on the immediately preceding state.  A consequence is that evaluation of prediction model \cref{eq:update-timemap} requires an iteration over time steps.  Thus the computational expense of this approach includes a linear scaling with $N_t$, in contrast to~\cite{HESTHAVEN201855}.

The change of model form from \cref{eq:timemap} to \cref{eq:update-timemap} also demands a change in the structure of the loss, though the data set defining the loss is essentially the same. %
For a $\mathcal{B} \subset \mathcal{F}$ we propose the loss
\begin{equation}\label{eq:loss_function}
  L_T(\mathcal{N}_T;\mathcal{B}) = \frac{1}{|\mathcal{B}|} \sum_{(k,\vmu)\in \mathcal{B}} \left\|\mathcal{N}_T(\vc(t_{k},\vmu),\vmu) - \vc(t_{k+1},\vmu)\right\|_2^2
\end{equation}
where the inputs to the model are the POD coefficients defined at $(t_k,\vmu)$ and the output is the coefficient defined at $(t_{k+1},\vmu)$.
Here again, the ideal map satisfies $\mathcal{N}_T =\argmin_{\mathcal{N}} L_T(\mathcal{N};\mathcal{F})$.
Comparing to the loss in~\cref{eq:timemap-loss}, minimizing this loss targets mappings learning the transition from a coefficient at time $t_k$ to $t_{k+1}$ for any parameter value.  Note that the deceptive simplicity of the loss is due to a judicious choice of training set notation designed to emphasize the serial nature of time evolution.

The loss~\Cref{eq:loss_function} measures the error by the network in taking a single timestep; however, the approximation of a time-dependent differential equation requires multiple time steps.  That is, the coefficients at time $t_k$ cannot simply be computed via $\mathcal{N}_T(\vc(t_{k-1},\vmu))\approx\vc(t_{k},\vmu)$, because $\vc(t_{k-1},\vmu)$ is not available for $k \geq 2$.  Instead, the coefficients must be computed from the initial condition, i.e.
\begin{align}
\begin{split}
	\vv_1(\vmu) &= \mathcal{N}_T(\vc(t_0,\vmu),\vmu)\\
	\vv_2(\vmu) &= \mathcal{N}_T(\vv_1(\vmu),\vmu)\\
	&\vdots\\
	\vv_{k-1}(\vmu) &= \mathcal{N}_T(\vv_{k-2}(\vmu),\vmu)\\
	\vc(t_k,\vmu) &\approx \mathcal{N}_T(\vv_{k-1}(\vmu),\vmu).
\end{split}
\end{align}
In fact, we require $\vv_j(\vmu) \approx \vc(t_j,\vmu)$ for $j = 1,2,\dots, N_t$.  To achieve this, we can adopt an alternative loss function to~\Cref{eq:loss_function} that takes multiple steps into account.

For any $(k,\vmu) \in \mathcal{B}$, the approximation coefficients at time $t_{k+m}$, where $m$ is a positive integer are computed by
\begin{subequations}
	\begin{align}
		&\vv_1^{(k)}(\vmu) \coloneqq \mathcal{N}_T(\vc(t_k,\vmu),\vmu),\\
		&\vv_p^{(k)}(\vmu) \coloneqq \mathcal{N}_T(\vv_{p-1}^{(k)}(\vmu),\vmu),\qquad p=2,\dots m,
	\end{align}
\end{subequations}
where $\vv^{(k)}_p$ denotes the approximation computed by $p$ evaluations of the neural network, starting with the exact value of the coefficients at time $t_k$.
For fixed $m\geq 1$, which serves as a hyperparameter of the training algorithm, we define the multi-step loss
\begin{equation}\label{eq:loss_function_multistep}
	L_T^m(\mathcal{N}_T;\mathcal{B}) = \frac{1}{|\mathcal{B}|} \sum_{(k,\vmu)\in \mathcal{B}} \sum_{p=1}^m\frac{1}{p}\|\vv_p^{(k)}(\vmu) - \vc(t_{k+p},\vmu)\|_2^2.
\end{equation}

\subsubsection{Residual Neural Networks}\label{sec:rnns}

As detailed in~\cite{carlberg2017galerkin}, the projected POD coefficients satisfy a low-dimensional dynamical system.  Thus, there is a solution operator $\mathcal{L}$ such that $\vc(t_{k+1};\vmu) = \vc(t_k;\vmu) + \int_{t_k}^{t_{k+1}}\mathcal{L}\left[\vc(t;\vmu)\right]\dif t$, so that $\|\vc(t_{k+1};\vmu)-\vc(t_k;\vmu)\| \leq C |t_{k+1}-t_k|$ for some constant $C$.  As a result, if the discrete time-steps are taken from a stable and accurate scheme, the difference $\vc(t_{k+1},\vmu) - \vc(t_k,\vmu)$ is relatively small.

\begin{figure}[!ht]
  \centering
  \includegraphics{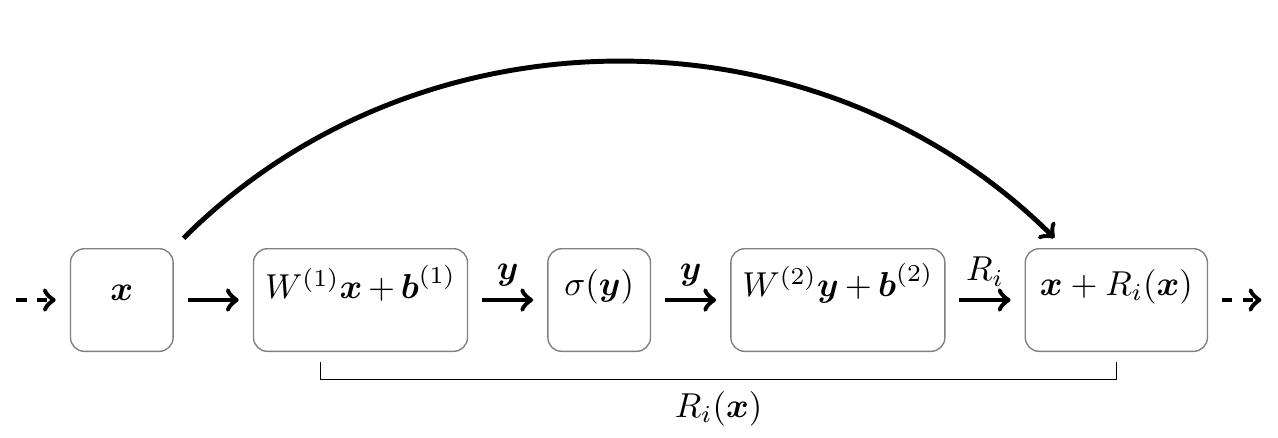}
  \caption{One residual layer of a ResNet showing the structure of the residual \textit{basic block} used.  The matrices and bias vectors defined the parameters for this layer.}\label{fig:res_block}
\end{figure}

Motivated by this perturbation, we consider the use of \textit{residual neural networks}~\cite{he2016deep} or ResNets based not only on their favorable training properties but also on their established connection to dynamical systems (see for instance~\cite{weinan2017proposal,haber2017stable,chen2018neural,wu2020data,patel2021physics}).  \Cref{fig:res_block} presents a schematic of single ResNet layer.  Each layer of the neural network computes a correction to the input, depicted by $R_i$ in the figure, combined with a \emph{skip} connection (the large curve in the figure).  To be concrete, for a layer $i$ with correction $R_i$ the ResNet architecture with $L$ layers is written as composition of identity mappings and correction terms
\begin{equation}
(I+R_{L-1}) \circ \ldots \circ (I+R_0) \circ x = x+R_0(x) + \ldots.
\end{equation}
The neural network weights are embedded in the $R_i$ mapping, which itself is expressed as a finite sequence of affine mappings and component-wise nonlinear transformations.  Specifically, if $R = R_i$ is the $i$th correction with $H$ hidden layers, then
\begin{equation}
	R(\vz) = A^{(H)} \circ \sigma \circ A^{(H-1)} \circ \sigma \circ \cdots \circ \sigma\circ A^{(1)}(\vz)
\end{equation}
where each $A^{(j)}$ is an affine transformation,
\begin{equation}
 A^{(j)} : \vz^{(j-1)} \mapsto W^{(j)}\vz^{(j)} + \vb^{(j)}
 \end{equation}
 and $\sigma$ is a component-wise mapping:
 \begin{equation}
 	\left[\sigma(\vz^{(j)})\right]_k = \sigma(z^{(j)}_k)
 \end{equation}

$W^{(j)}$ denotes the weight matrix for layer $j$, $\vb^{(j)}$ denotes the bias vector for layer $j$; these matrices and vectors are the parameters in the neural network to be determined by training.  The dimensions of the weight matrices and biases are chosen so that each composition is well-defined, and the final output dimension must match the input dimension, so that the sum $\vz + R(\vz)$ can be computed.
We point out that the choice of neural network architecture within $R_i$ can take a more general form.  Other examples for image classification can be found in~\cite{he2016deep}.

\begin{figure}[!ht]
\centering
  \includegraphics[width=\textwidth]{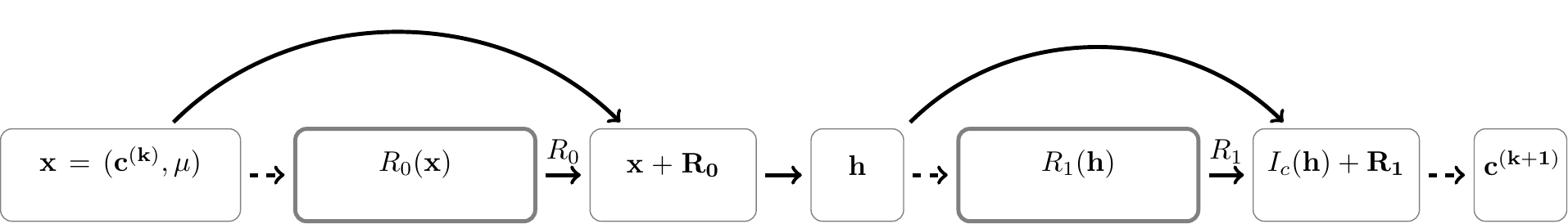}
  \caption{A two layer ResNet for time-stepping (see also~\cref{eq:resnet-timestep-update}).  The residual $R_i$ functions contain the weights and biases~--~cf.~\cref{fig:res_block}.}\label{fig:initial_res_block}
\end{figure}

For our time-stepping approach the input to the ResNet is defined as $\vx = (\vc,\vmu)\in \mathbb{R}^{N_{\rb} + P}$.  In a modest modification of the
standard ResNet, where outputs and inputs are the same size, in the final layer of our ResNet we include a contraction operator $I_c(\vc,\vmu) \rightarrow \vc$,
which ensures the output is equal to the size of the coefficient vectors.  Including the contraction operator with, the coefficient mapping described by~\cref{eq:update-timemap} is defined as
\begin{equation}\label{eq:resnet-timestep-update}
\mathcal{N}_T(\vc,\vmu) := (I_c+R_{L-1}) \circ \ldots \circ (I+R_0) \circ (\vc,\vmu)
\end{equation}
where $I_c(\cdot)\in \mathbb{R}^{N}$ and $R_{L-1}(\cdot) \in\mathbb{R}^{N}$ provide the necessary contraction. \Cref{fig:initial_res_block} is a graphical representation of this network for $L=2$.

\subsection{Training the Residual Network}
To train the residual network $\mathcal{N}_T$, we select a parameter set $\mathcal{D}_{\text{train}} = \{\vmu_1,\vmu_2,\dots, \vmu_{N_{\text{train}}}\}$ and construct two collections of features/targets:
\begin{subequations}
\begin{align}
  \mathcal{F}_{\text{train}} & = \bigg\{\left(k, \vmu_i\right) : k = 0,\dots,N_t-1,\ \ i=1,\dots,N_{\text{train}}\bigg\}\\
  \mathcal{T}_{\text{train}} & = \bigg\{ \vc(t_{k+1};\vmu_i) : k = 0,\dots,N_t-1,\ \ i=1,\dots,N_{\text{train}}\bigg\}
\end{align}
\end{subequations}
For a fixed $m\geq 1$, the loss function $L_T^m(\mathcal{N}_T;\mathcal{F}_{\text{train}})$ from~\Cref{eq:loss_function_multistep} is a measure of the difference between the observed coefficients and the predicted coefficients returned by the network.

The optimal network $\mathcal{N}_T$ is found by minimizing this loss.  We use the L-BFGS optimization method~\cite{liu1989limited, nocedal1980updating} to determine optimal weight and biases in our examples, but other common optimization routines, including variants of the stochastic gradient descent algorithm, can be used as well.  In the latter case, each step produces a descent direction for a randomly selected batch $\mathcal{B} \subset \mathcal{F}_{\text{train}}$ using back-propagation to compute the gradient that increments the weights.  For L-BFGS, much larger batch sizes must be used; we take $\mathcal{B} = \mathcal{F}_{\text{train}}$.

Since the POD coefficients are constructed to be orthonormal in the $L^2$-inner product, this corresponds to the $L^2$-norm of the difference of the projection-based solution and the neural network solution~---~if
$u^{\text{rb}}$ the projection-based RB solution and $u^{\mathcal{N}_R}$ represents the solution constructed from the neural network mapping, then~\cref{eq:loss_function} or~\cref{eq:loss_function_multistep} computes
$\|u^{\text{rb}} - u^{\mathcal{N}_R}\|_{L^2}$.

If $m \ll N_t$, then the training set can also be used as a validation set, e.g, when using multiple restarts or early stopping.  For every $\vmu$ in the training set, the neural network is evaluated $N_t$ times on the initial POD coefficients; this is our approach in numerical examples.

\section{Temporal Reduced Basis Algorithm}\label{sec:algorithm}

We next summarize the temporal reduced basis algorithm introduced over the previous sections.
The offline stage of reduced order modeling is split into two parts.  First the traditional POD basis construction, summarized in~\cref{alg:pod},
takes a matrix of FE degrees of freedom, denoted $A$, and tolerance $\epsilon$ as inputs and computes a POD with basis functions that are $L^2$-orthonormal~---~if $M_{ij} = (\phi_j,\phi_i)_{L^2}$ is the mass matrix, the output $V = \text{POD}(A,\epsilon)$ satisfies $V^T M V = I$.  For a detailed description of the traditional POD algorithm, see Section 6.3 of~\cite{quarteroni2015reduced}.
\begin{algorithm2e}[!ht]
  \DontPrintSemicolon%
  \SetKwProg{Fn}{}{}{end}
  \SetKwFunction{POD}{POD}
  \SetKwFunction{SVD}{SVD}

  \BlankLine%
  $M_{ij} = (\phi_j, \phi_i)_{L^2}$ for $i,\, j=1,\dots N_h$\tcc*{Mass matrix}
  \Fn{\POD{$A$, $\epsilon$}}{%
    $U \Sigma V^T \leftarrow \SVD(M^{1/2}A)$\tcc*{Compute the SVD of $N_h \times m$ matrix $M^{1/2}A$}
    $\sigma_k = \left(\mathrm{diag}(\Sigma)\right)_k$\;
    $\hat{m} = $ smallest $m$ such that $\frac{\sum_{k=m+1}^{r}\sigma_k^2}{\sum_k \sigma_k^2} < \epsilon$\tcc*{Find $\hat{m}$ to a tolerance}
    $\widehat{U} = M^{-1/2}U$\tcc*{Transform to $L^2$-orthonormal basis}
    $\hat{A} = \begin{bmatrix}\widehat{U}_{*1} & \widehat{U}_{*2} & \cdots & \widehat{U}_{*\hat{m}}\end{bmatrix}$\tcc*{Extract $\hat{m}$ POD basis vectors}
    \Return{$\hat{A}$}
  }
  \caption{POD basis construction.}\label{alg:pod}
\end{algorithm2e}

The traditional POD basis construction is used inside the algorithm that constructs the training data, summarized in~\cref{alg:trainingdata}.  The parameter space $\mathcal{D}$ is sampled to produce a finite set of parameters, $\mathcal{D}_\train$.  In this paper, the training set was built by sampling from a uniform grid; for high-dimensional parameter domains, a uniform grid may be computationally intractable, and the training sample should be constructed using a random sample.\par

The full order solution is computed for each parameter in the training and validation set, and the initial condition is subtracted from the vector of degrees of freedom.  Thus, the neural network will learn how the POD coefficients of variation from the initial condition propagate in time.  For this paper, we do not store every time step computed during the simulation.  This is done to emulate cases where only a subset of time steps are saved; for example, when experimental data is available instead of a numerical model, or using adaptive time-stepping techniques.\par

We use the training data to construct a reduced basis using the two-step POD process described in~\cref{sec:pod_construction}.  Thus, we perform a compression in time, followed by a compression in parameter space.  We finally store the true POD coefficients of the projected full-order solutions for all parameters in the training and validation set.
\begin{algorithm2e}[!ht]
  \SetArgSty{textrm}
  \SetKwFunction{LHS}{LatinHypercube}
  \SetKwFunction{grid}{grid}
  \SetKwFunction{POD}{POD}
  \DontPrintSemicolon%
  \BlankLine%
  $\mathcal{D} \leftarrow \mathcal{D}_1 \times \cdots \times \mathcal{D}_P$\tcc*{Decomposition of parameter space}
  $\Ntrain \leftarrow N_1 * \cdots * N_P$\tcc*{Number of parameter samples for each parameter}
  $\Dtrain \leftarrow \grid(\mathcal{D},\, \Ntrain)$\tcc*{Full \grid\ or \LHS\ sample of $\Ntrain$ points of $\mathcal{D}$}
  \;
  \ForEach(\tcc*[f]{Compute the full order solution for each parameter}){$\vmu_i \in \mathcal{D}_\train$}{
    \For(\tcc*[f]{\dots\ at each time step}){$t_k = k \Delta t =0\dots T$}{
      \If(\tcc*[f]{Save FE coefficients every $\bar{k}$ time steps}){$k\bmod \bar{k} = 0$}{
        $\hat{k} = \frac{k}{\bar{k}}$\tcc*{Compute the index for the $k^{\text{th}}$ step}
      $\left(U_\train\right)_{i\hat{k}j} \leftarrow \alpha_j(t_k;\, \vmu_i)
      -\alpha_j(0;\, \vmu_i)
      $
      \quad for $j=1,\dots,N_h$\tcc*{Scaled FE coefficients}
    }
    }
  }
  \;
  \ForEach(\tcc*[f]{Compress the solution in time with POD}){$\vmu_i \in \mathcal{D}_\train$}{
    $T_i = \left(U_\train\right)_{i,*,*}^T$\tcc*{Extract the $N_h \times m=\frac{N_h}{\bar{k}}$ $i^\text{th}$ submatrix}
    $\hat{T}_i = \POD(T_i, \epsilon_t)$\tcc*{Construct a reduced $N_h\times \hat{m}$ POD basis}
  }
  $T = \begin{bmatrix} T_1 & T_2 & \cdots & T_P\end{bmatrix}$\tcc*{Agglomerate $T_i$ over parameters}
  $W = \POD(T, \epsilon_{\vmu})$\tcc*{Compress over parameters with POD}
  \;
  $\vc(t_k;\, \vmu_i^{(\train)}) =  W^T M \left(U_\train\right)_{i,\hat{k},*}$\tcc*{Compute the POD coefficients via $L^2$ projection}\;
  \caption{Construct training data.}\label{alg:trainingdata}
\end{algorithm2e}

The final step of the offline stage is to train the neural network; this is presented in~\cref{alg:training}.  The number of steps $m$ is treated as a hyperparameter of the algorithm, which determines the number of consecutive network evaluations that contribute to the loss.  For every $(k,\vmu)$ pair in a training batch, the network is evaluated $m$ times to obtain an approximation for $\vc(t_{k+m},\vmu)$.  The weighted sum of~\cref{eq:loss_function_multistep} is used to capture the error in the coefficients for $t_{k+1},\dots,t_{k+m}$.

Algorithms for training neural networks are often sensitive to the initialization of weights and biases, both on training loss and generalization to new data~\cite{goodfellow2016deep}, so it is an effective strategy to train multiple networks with different random initializations, and choose the one that performs best on a validation set.  This is known as a multiple restarts approach~\cite{HESTHAVEN201855, WaHeRa_2019_nonintrusivenn}.~\cref{alg:training} takes this approach, and selects the network with lowest average validation error.

The training set is not usually appropriate to use for model validation, because the network has been exposed to this data.  However, as long as $m \ll N_t$, a propagation of the coefficients for the full number of time steps is sufficient to ensure the network cannot simply memorize the data.
\begin{algorithm2e}[!ht]
  \SetArgSty{textrm}
  \SetKwFunction{LHS}{LatinHypercube}
  \SetKwFunction{grid}{grid}
  \SetKwFunction{POD}{POD}
  \DontPrintSemicolon%
  \BlankLine%
  \Input{$\Dtrain$, $\vc(t_{\hat{k}};\, \vmu_i^{(\train)})$}
  \For{$\ell = 1,2,\dots, \text{restarts}$}{

  Initialize neural network $\mathcal{N}_T^{(\ell)}$

  \For{$j=1,\dots,\text{epochs}$}{\For{batch $\mathcal{B} \subset \mathcal{F}_{\text{train}}$}{\For{$(k,\vmu)\in\mathcal{B}$}{$\vv_1^{(k)}(\vmu) = \mathcal{N}_T^{(\ell)}(\vc(t_k,\vmu),\vmu )$

  $\vv_p^{(k)}(\vmu) = \mathcal{N}_T^{(\ell)}(\vv_{p-1}^{(k)}(\vmu),\vmu )$\quad for $p=2,\dots, m$ \tcc*{Propagate coefficients with network}}Evaluate $L_T^m(\mathcal{N}_T^{(\ell)};\mathcal{B}) = \frac{1}{\vert\mathcal{B}\vert}\sum\limits_{k,\mu\in \mathcal{B}}\sum\limits_{p=1}^m\frac{1}{p}\left\|\vv_p^{(k)}(\vmu) - \vc(t_{k+p},\vmu)\right\|_{\mathbb{R}^{N_{rb}}}^2$ \tcc*{See~\cref{eq:loss_function_multistep}}

  Adjust weights and biases via optimization step}}

  \For{$\vmu \in \mathcal{D}_\train$}{$\vv_1(\vmu) = \mathcal{N}_T^{(\ell)}(\vc(t_0,\vmu),\vmu )$

  $\vv_p(\vmu) = \mathcal{N}_T^{(\ell)}(\vv_{p-1}(\vmu),\vmu )$\quad for $p = 1,\dots, N_t$\tcc*{Propagate coefficients until time $T$}}

  $\hat{\mathcal{C}}_\ell \leftarrow \frac{1}{\Ntrain}\sum\limits_{\vmu\in \mathcal{D}_\train}\left \| \vv_{N_t}(\vmu) - \vc(t_{N_t};\vmu_i) \right\|^2_{\mathbb{R}^{N_{rb}}}$\tcc*{Compute average error at final timestep}}
  $\ell^* = \argmin_\ell \hat{C}_\ell$

  \Output{$\mathcal{N}_T^{(\ell^*)}$}
  \caption{Train the neural network.}\label{alg:training}
\end{algorithm2e}

\section{Numerical Evidence}\label{sec:numerics}

We apply our proposed method to obtain ROMs for several time-dependent parametrized PDEs. For each PDE, we select a set of training and test parameters, resulting in well-resolved, full-order FEM solutions for all parameters. For all simulations, we use piecewise linear elements and the Firedrake~\cite{firedrake} library to construct the discretization.

\subsection{Advection-Diffusion Model}

The first two examples present a proof of concept demonstration that the proposed neural network prediction is effective in the online stage.
To this end, both problems admit an affine decomposition, so an efficient offline-online procedure may be attained using standard Galerkin projection RB methods.

\subsubsection{One-Dimensional Case}

Following the example in~\cite{parish2020time}, we consider the parameterized linear advection-diffusion equation given as
\begin{subequations}\label{eq:example-1}
\begin{align}
  u_t + \mu_1u_x &= \mu_2 u_{xx}, \qquad x\in (0,2),\ t \in (0,T]\\
  u(0,t) &= 0 = u(2,t)\\
  u(x,0) &= x(2-x)e^{2x}
	\end{align}
\end{subequations}
with final time $T = 0.3$, advection speed $\mu_1 \in [-2, -0.1]$, and diffusion coefficient $\mu_2 \in [0.1, 1]$.
The full-order simulation is spatially discretized on a uniform mesh with $h=2/101$.
Temporally, the discretization uses Crank-Nicholson time-stepping with  $\Delta t = 3\times 10^{-4}$ for $N_t = 1000$ steps until the final time $T$.

The training set is generated using a $9\times 9$ uniform grid of $(\mu_1, \log_{10}(\mu_2)) \in [-2, -0.1] \times [-1, 0]$, leading to $N_{\text{train}} = 81$ samples of the parameter domain $\mathcal{D} = [-2, -0.1]\times [0.1, 1.0]$.  To imitate the situation where a high-fidelity solution is only saved to disk periodically, we use the data every 10 timesteps~---~that is, we only have access to the solution at times $t_0$, $t_{10}$, $t_{20}$, $\dots$, $t_{990}$, $t_{1000}$.

Thus, the collection of features is
\begin{equation}
	\mathcal{F}_{\text{train}} = \left\{(\hat{k},\vmu_i): \hat{k}=0,\dots, 99,\ i = 1,\dots, 81  \right\},
\end{equation}
Where $\vmu_i \in \mathcal{D} = [-2, -0.1]\times [0.1, 1]$, $\hat{k} = \frac{k}{10}$, and for each $\hat{k}$ and $\vmu_i$, there is a corresponding full-order solution
\begin{equation}
	\valpha(t_{10\hat{k}};\vmu_i) = \valpha(t_k;\vmu_i),
\end{equation}
where $\valpha(t_k;\vmu_i)$ is a vector of degrees of freedom representing the solution $u(x,t_k;\vmu_i)$ to~\cref{eq:example-1}.

The reduced basis is constructed by performing the POD process in two stages~\cite{WaHeRa_2019_nonintrusivenn} on the training parameters.  As in~\cite{parish2020time}, the proper orthogonal decomposition is constructed with respect to the deviation from the initial condition: i.e., degrees of freedom are stored in an array  $(U_{\train})_{i\hat{k}j} = \alpha_j(t_k;\vmu_i) - \alpha_j(0;\vmu_i)$, and the two-step process outlined in~\cref{alg:trainingdata} is performed.

We choose the tolerance parameters to be $\epsilon_t = 10^{-4}$ (compression in time) and $\epsilon_{\vmu} = 10^{-4}$ (compression over parameters), resulting in 8 $L^2$-orthonormal basis functions.  The $L^2$-projection coefficients are computed for every parameter pair in the training set for $t_0$, $t_{10}$, $t_{20}$, $\dots$, $t_{990}$, $t_{1000}$.  The corresponding coefficient targets are computed as detailed in~\cref{alg:trainingdata}, resulting in the collection of targets
\begin{equation}
	\mathcal{T} = \left\{\vc(t_{\hat{k} + 1}; \vmu_i) : \hat{k} = 0,\dots, 99,\ i=1,\dots, 81  \right\}.
\end{equation}

We train a time-stepping network with two residual blocks, followed by a linear output layer using~\cref{alg:training}.  The first block has input dimension 10 (eight coefficients and two parameters), with output dimension 8.  The block has two hidden layers of width $W = 10$.  The second block has input and output dimension of 8; the hidden and output layers have width $W = 8$.  We use the exponential linear unit~\cite{clevert2015fast} as the nonlinear activation function, and optimize using the L-BFGS method~\cite{nocedal1980updating}.  A full listing of parameters for the network architecture can be found in the Appendix.

In the results that follow, the quality of the online prediction is measured on
 a test set of $N_{\text{test}} = 50$ parameters in $[-2,-0.1]\times [0.1,1]$, generated using Latin hypercube sampling.  For each parameter, we compute the full-order solution of~\cref{eq:example-1}, and the corresponding projection coefficients for $t_{10} < t_{20}< \cdots< t_{990}< t_{1000} = 0.3$.  We then compute approximate coefficients using the trained neural network, along with the coefficients using the classical POD method, following the same discrete time-steps as the projection.

\begin{figure}[!ht]
  \centering
  \includegraphics{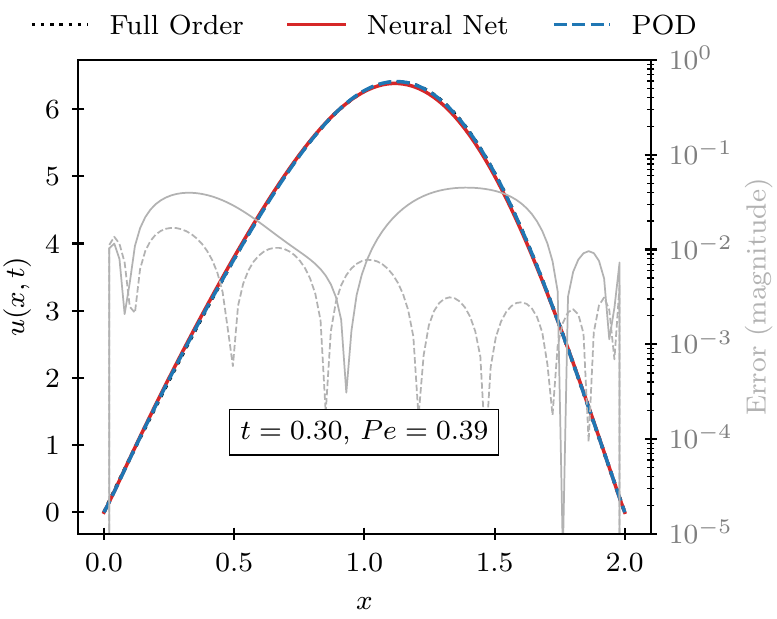}
  \hfill
  \includegraphics{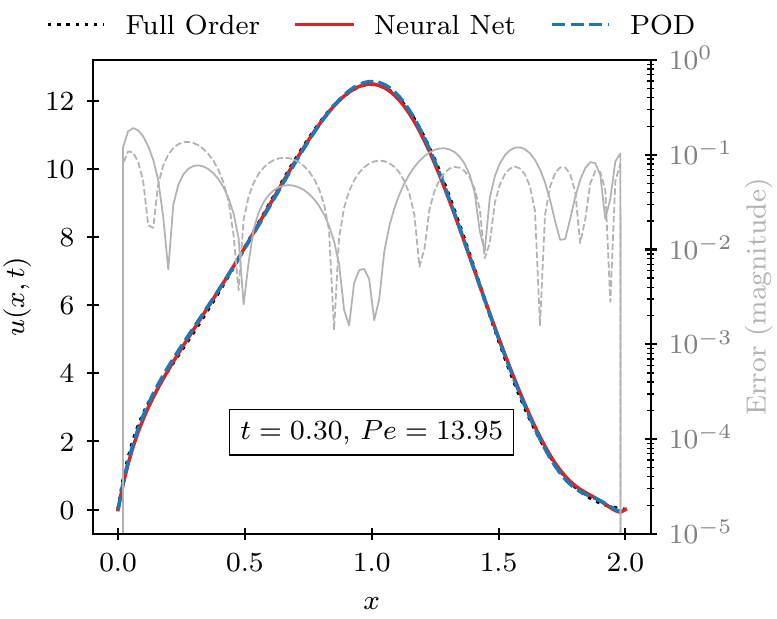}
  \caption{\textbf{1D advection-diffusion}: Two solution profiles at different parameter values.  The pointwise error (in magnitude) with the full-order solution in the neural network approximation (solid) and the POD (dashed) is shown in the background.}\label{fig:advec-diff-1d-profile}
\end{figure}

\Cref{fig:advec-diff-1d-profile} shows solutions of the PDE at two distinct points in parameter space: the full-order solution, the neural network prediction, and the classical POD approach.  The two points correspond to Peclet numbers of $Pe=0.39$ and $Pe=13.95$, which are the smallest and largest in the parameter test set.
The solution profile shows that there is variation between the solutions, with the high Peclet number case corresponding to a less diffusive solution.  Despite this variability the ROM with a neural network reconstruction (red) is able to produce an accurate approximation of the solution profile in both cases.  In the case of $Pe = 0.39$, the classical POD solution attains higher accuracy over most of the computational domain.  However, with more dominant advection ($Pe = 13.95$), the accuracy between the two methods is more comparable, with the neural network solution providing a more accurate approximation over a significant part of the domain.

\begin{figure}[!ht]
  \centering
  \includegraphics{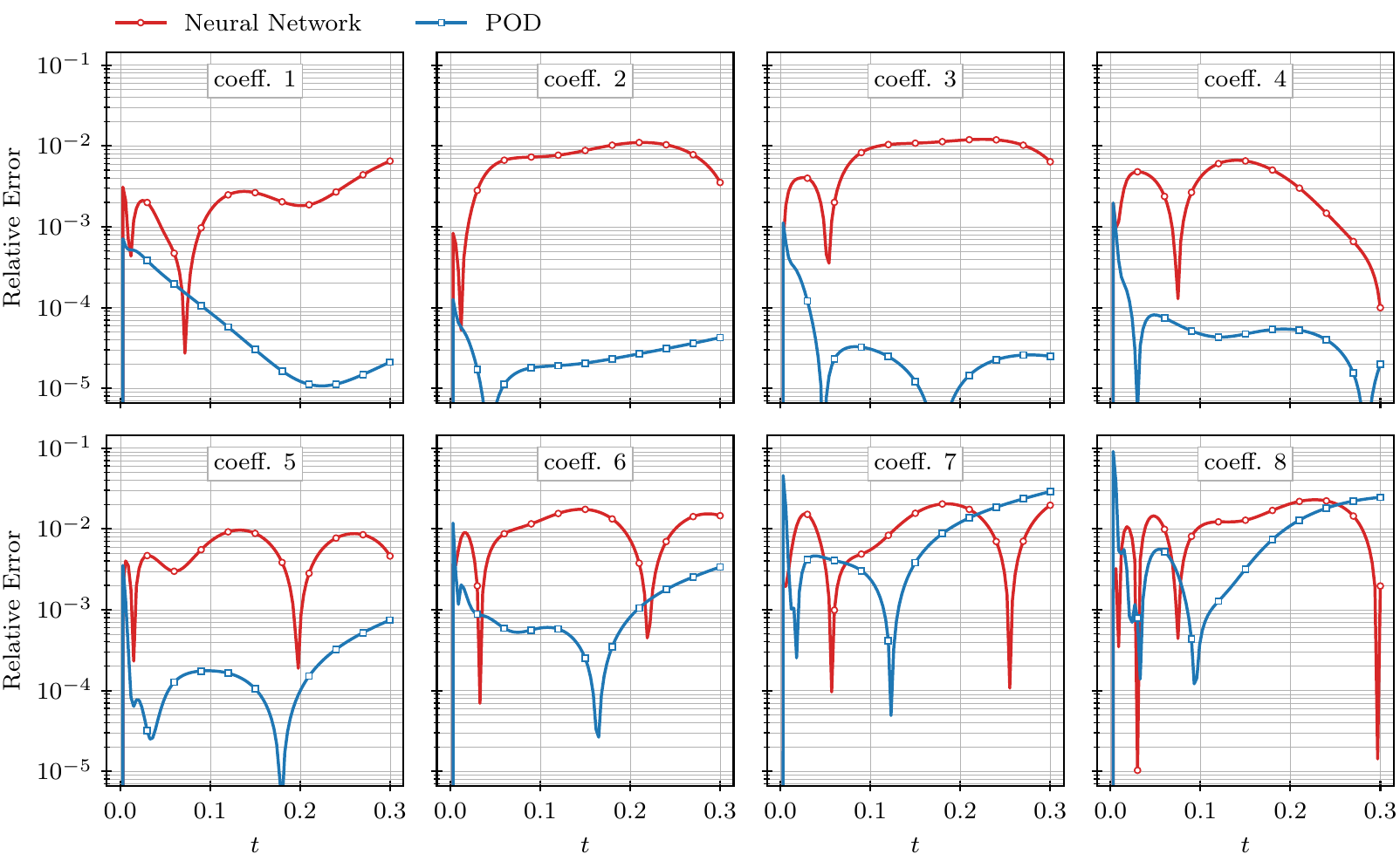}
  \caption{\textbf{1D advection-diffusion}: ($Pe=0.39$) The relative error with respect to the $L^2$-projection coefficients.  The red curve is the error in the neural network approximation and the blue curve shows the error in the coefficients computed using Galerkin POD\@.  Each plot corresponds to the error in one of the 8 POD coefficients, as indicated. The timesteps are marked every 10 steps.}\label{fig:ex_1_coeff_err_min_peclet}
\end{figure}

\begin{figure}[!ht]
  \centering
  \includegraphics{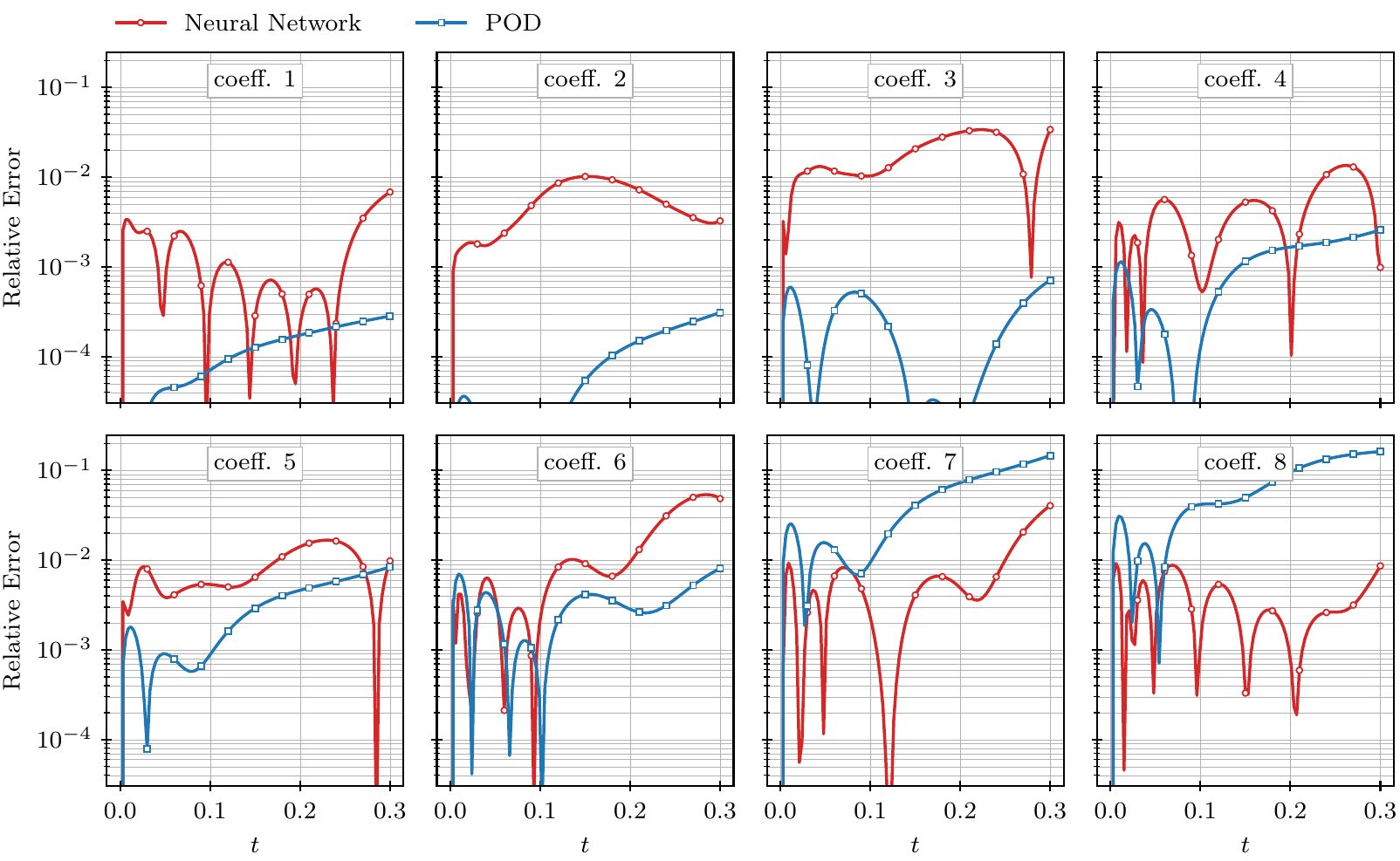}
  \caption{\textbf{1D advection-diffusion}: ($Pe=13.95$) The relative error with respect to the $L^2$-projection coefficients.  The red curve is the error in the neural network approximation and the blue curve shows the error in the coefficients computed using Galerkin POD\@.  Each plot corresponds to the error in one of the 8 POD coefficients, as indicated.  The timesteps are marked every 10 steps.}\label{fig:ex_1_coeff_err_max_peclet}
\end{figure}

\Cref{fig:ex_1_coeff_err_min_peclet} and~\cref{fig:ex_1_coeff_err_max_peclet} show
the accuracy in the ROM approximation measured by the error of the coefficients with respect to the $L^2$-projection of the full-order solution. %
This shows the error in the projection coefficients for $Pe = 0.39$ and $Pe = 13.95$, respectively.
For the diffusion-dominated case~\cref{fig:ex_1_coeff_err_min_peclet}, the coefficients computed via the POD Galerkin method are a better approximation to the true projection coefficients, except for the last two near the end of the time interval in question.
In contrast,~\cref{fig:ex_1_coeff_err_max_peclet} the approximation error in the neural network is much closer to the error in the POD Galerkin solution, and captures the higher frequency coefficients with higher accuracy.

\begin{figure}[!ht]
  \centering
  \includegraphics{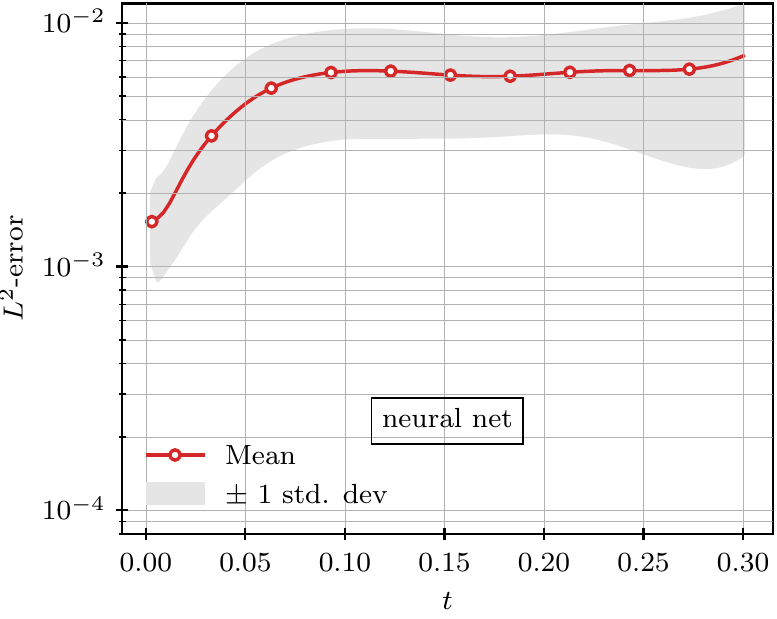}
  \hfill
  \includegraphics{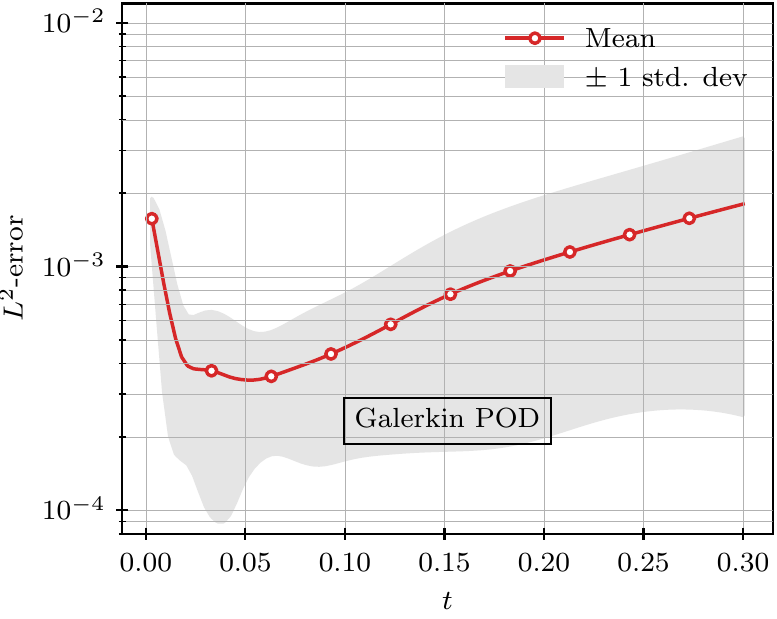}
  \caption{\textbf{1D advection-diffusion}: (Left) Mean relative $L^2$-error with respect to full-order solution. The shaded region denotes $\pm$ one standard deviation over the test set.  (Right) Relative $L^2$-error for the parameter with the largest mean error. In each case, timesteps are marked every 10 steps.}\label{fig:ex_1_mean_err}
\end{figure}

\Cref{fig:ex_1_mean_err} explores the accuracy of the neural network prediction over the entire test set as a function of time.
The image on the left plots the mean, and a single standard deviation envelope of the relative $L^2$-error. This shows that the test error
remains bounded below $4\%$ error for the majority of the samples in the test space. In fact, at the final 
time $T = 0.3$, which corresponds to 100 timesteps of size $10\Delta t$, the mean relative $L^2$ error with respect to the true projection is $7.3\times 10^{-3}$.  Compared with POD Galerkin, the neural network has higher relative error over the test parameters on average.  However, the error in the POD Galerkin method increases over time, while the neural network error is more nearly constant over the time interval.

\subsubsection{Two-Dimensional Case}

Our second example considers the two-dimensional advection-diffusion problem
\begin{subequations}\label{eq:example-2}
\begin{align}
	u_t - a\Delta u + \vb(\vmu)\cdot\nabla u = 0&, \quad \vx \in \Omega,\ t \in (0,T]\\
	\nabla u\cdot\vn = 0&, \quad \vx \in \partial\Omega \\
	u(\vx,0) = \exp(-10(x^2 + y^2))&.
\end{align}
\end{subequations}
Here $\vx = (x,y) \in \Omega = (0,1)^2$, the vector field $\vb(\vmu)$ is given by
\begin{equation}
\vb(\vmu) = \vb(\mu) = b\cdot(\cos(\mu),\sin(\mu))^T,
\end{equation}
where $\mu \in [0, \pi/2]$, $T = 0.5$, and the constants $a$ and $b$ are fixed at $0.5$ and $2$, respectively.  The full-order simulation is spatially discretized on a uniform mesh with $2,048$ linear triangular elements.  Temporally, we use backward Euler time-stepping with $\Delta t = 10^{-3}$ for $N_t = 500$ time steps.

The training set is composed of $N_{\train} = 50$ uniformly spaced values in the interval $[0, \frac{\pi}{2}]$, and again, only the data corresponding to $t_{0},t_{10}, t_{20},\dots, t_{490}, t_{500}$ are used.
The two-stage POD process is performed on $\valpha(t_k;\vmu) - \valpha(0;\vmu)$ with $\epsilon_t = \epsilon_{\mu} = 10^{-4}$, and results in 7 basis functions.

\begin{figure}[!ht]
  \centering
  \includegraphics{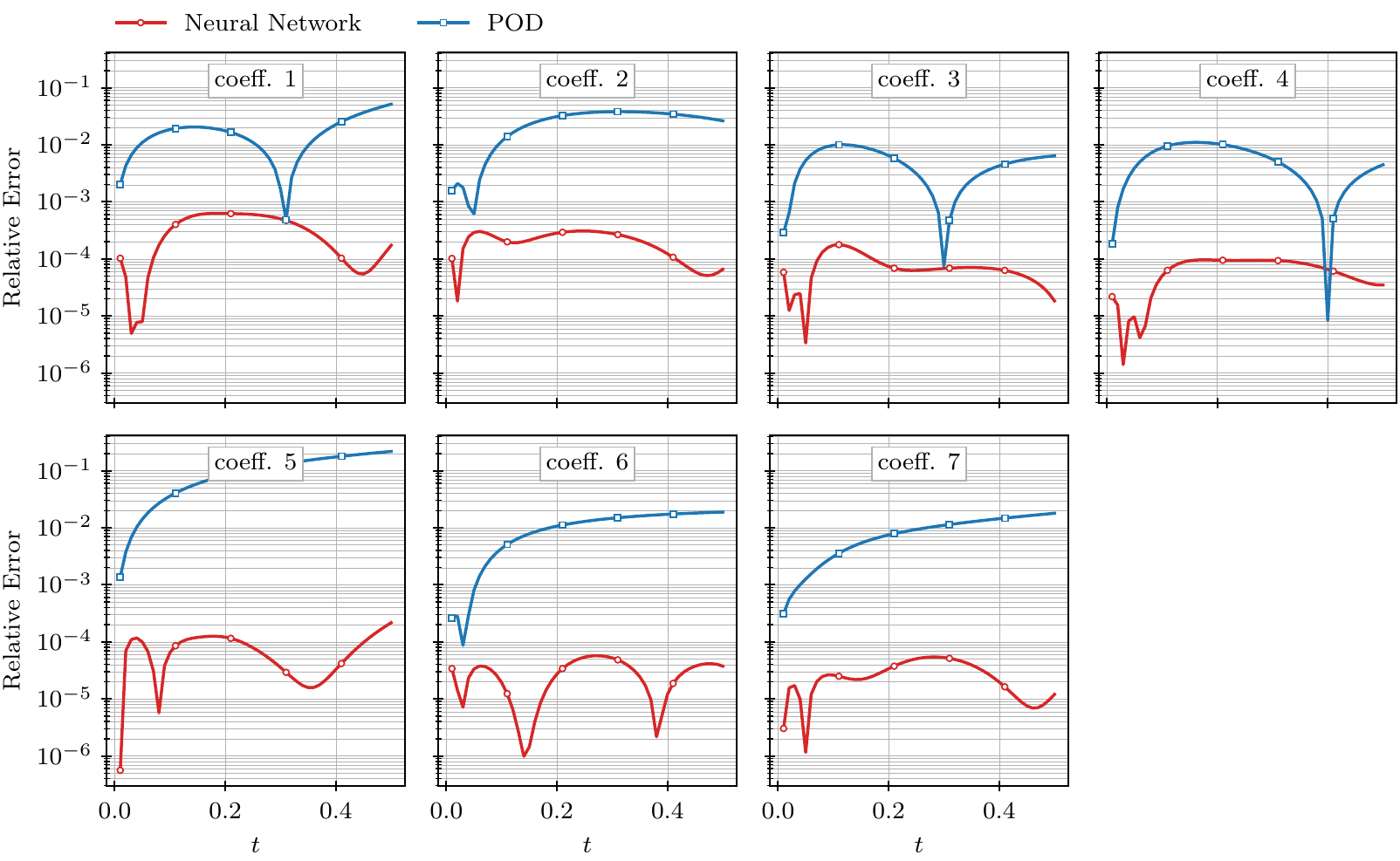}
\caption{\textbf{2D advection-diffusion}: The relative error with respect to the $L^2$-projection coefficients.  The red curve is the error in the neural network approximation and the blue curve shows the error in the coefficients computed using Galerkin POD\@.  Each plot corresponds to the error in one of the 7 POD coefficients, as indicated.  The timesteps are marked every 10 steps.}\label{fig:ex_2_coeff_err}
\end{figure}

For this example, the neural network solution results in much higher accuracy than the POD Galerkin projection.  In~\Cref{fig:ex_2_coeff_err} the absolute error with respect to the true projection coefficients is shown for $\mu \approx 0.49\times \pi$, which is the test parameter with best performance for the POD method.  The neural network consistently provides a more accurate approximation of the projection coefficients than POD\@.

\begin{figure}[!ht]
  \centering
  \includegraphics{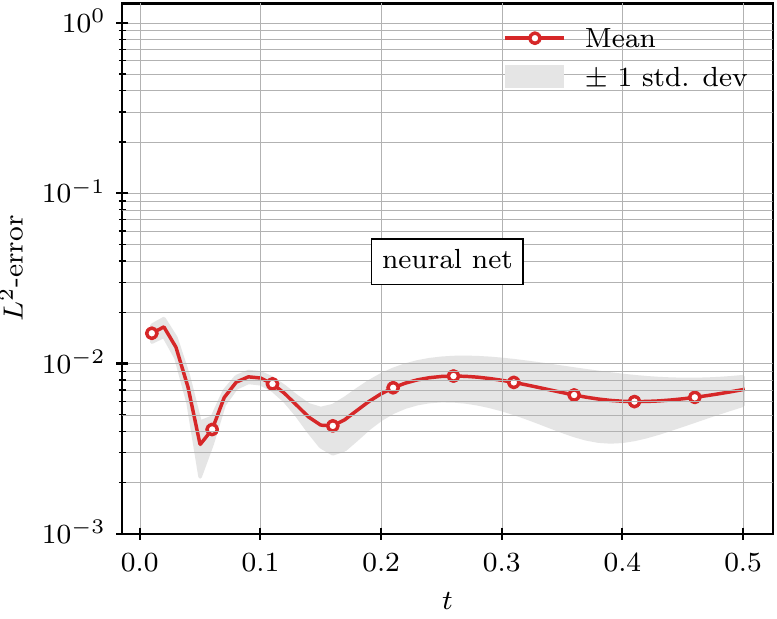}
  \hfill
  \includegraphics{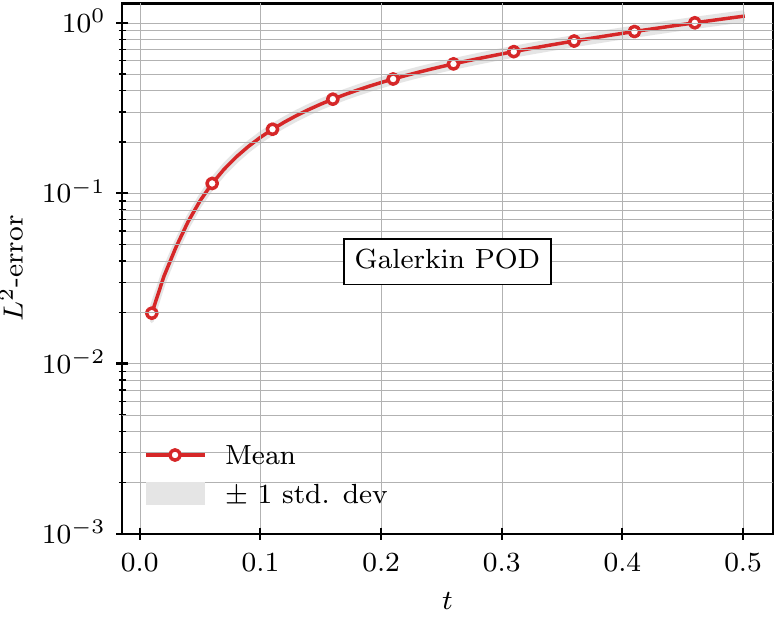}
  \caption{\textbf{2D advection-diffusion}: (Left) Mean relative $L^2$-error with respect to full-order solution. The shaded region denotes $\pm$ one standard deviation over the test set.  (Right) Relative $L^2$-error for the parameter with the largest mean error. In each case, timesteps are marked every 10 steps.}\label{fig:ex_2_mean_err}
\end{figure}

In~\Cref{fig:ex_2_mean_err}, the mean $L^2$-error over the testing set is plotted as function of time.  Galerkin POD is unable to capture the dynamics of the projection coefficients for the entire testing set.  In fact, classical POD techniques struggle to resolve advection-dominated problems, and must be augmented by various stabilization techniques in order to be effective~\cite{quarteroni2015reduced, giere2015supg, ohlberger2015reduced, rubino2018streamline}. On the other hand, the neural network is able to maintain accuracy over the time interval.  This demonstrates that our proposed method can be effective for problems where classical POD techniques fail.

\subsection{Non-affine Linear Diffusion-Reaction Equation}

\begin{figure}[!ht]
  \centering
  \includegraphics{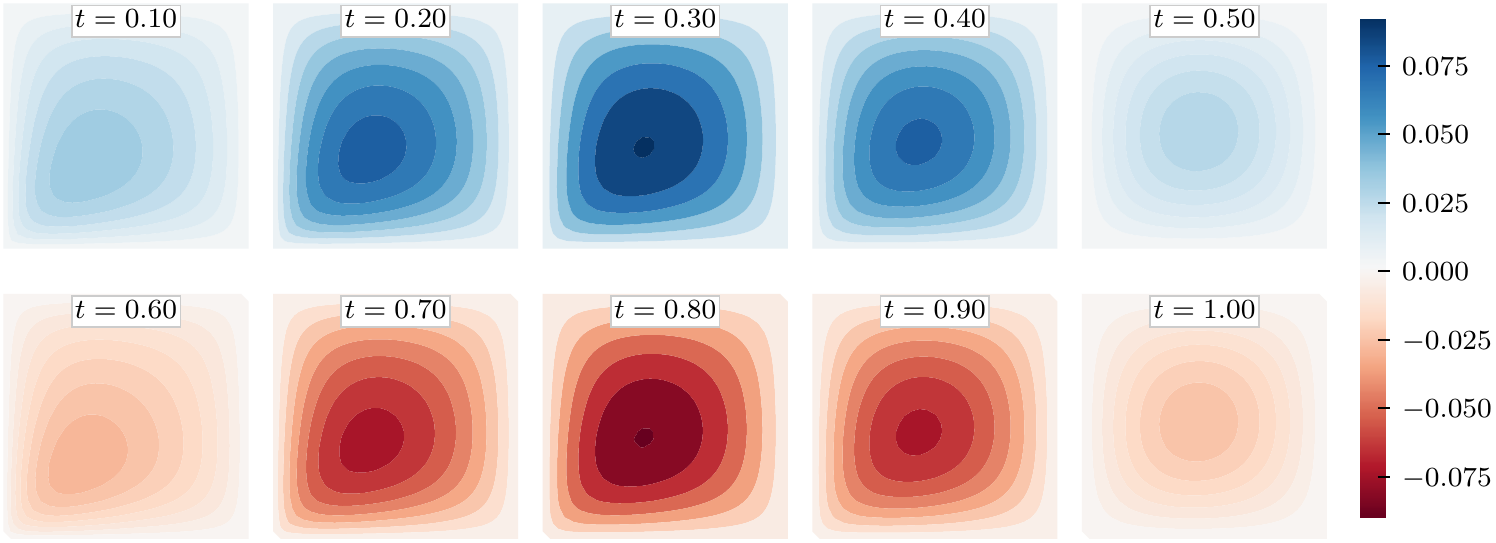}
  \caption{\textbf{Non-affine diffusion-reaction}: Solution snapshots at various points in time.}\label{fig:non-affine-solution-allsteps}
\end{figure}

We next apply our method to a two-dimensional diffusion-reaction equation that does not admit an affine decomposition~\cite{grepl2012certified, grepl2007efficient}:
\begin{subequations}\label{eq:nonaffine}
\begin{align}
	u_t - \Delta u +  g(\vx;\vmu)u = \sin(2\pi t)g(\vx;\vmu)&, \quad \vx \in \Omega,\ t \in (0,T]\\
	u = 0&, \quad \vx \in \partial\Omega \\
	u(\vx,0) = 0&, \quad \vx \in \Omega.
\end{align}
\end{subequations}
Here $\vx = (x,y)$, $\Omega = (0,1)^2$, and the function $g$ has the form:
\begin{equation}
	g(\vx;\vmu) = g(\vx;\mu_1,\mu_2) = \frac{1}{\sqrt{(x - \mu_1)^2 + (y - \mu_2)^2}},
\end{equation}
for $\vmu = (\mu_1, \mu_2) \in [-1, -0.01]^2$.
As an example, consider the solution to~\cref{eq:nonaffine} in the case of
$(\mu_1,\mu_2)=(-0.05,-0.05)$.  The solution oscillates in time, as depicted
in~\cref{fig:non-affine-solution-allsteps}.

The weak form of this problem is given by: find $u \in H_0^1(\Omega)$ such that
\begin{equation}
	(u_t, v)_0 + (\nabla u, \nabla v)_0 + (g(\vx;\vmu)u,v)_0 = \sin(2\pi t)(g(\vx;\mu),v)_0, \quad \forall v\in H_0^1(\Omega).
\end{equation}
The bilinear and linear forms
\begin{align}
	\begin{split}
		(g(\vx;\vmu)u,v)_0 &= \int_{\Omega} g(\vx;\vmu)uv\ \dif \vx,\\
		(g(\vx;\mu),v)_0 &= \int_{\Omega} g(\vx;\vmu)v\ \dif \vx,
	\end{split}
\end{align}
do not admit an affine decomposition.  To apply the Galerkin-POD method, the corresponding discrete operators must be assembled for every new parameter instance, which is prohibitively expensive for the online stage.  Alternative methods to recover approximate affine decompositions, such as the Empirical Interpolation Method, have been proposed.  However, these techniques require an additional offline phase where the parameter space is sampled, and their implementation is highly dependent on the problem at hand, and may require modifications of the underlying computational code~\cite{chaturantabut2010nonlinear,drohmann2012reduced, HESTHAVEN201855, casenave2015nonintrusive}.

For the training set, we generate a $10\times 10$ uniform grid of the interval $[-2,0]^2$, to generate samples $\vlambda = (\lambda_1, \lambda_2)$, and set $\mu_i = -10^{\lambda_i}\in [-1, -0.01]$.  In this case, $N_{\text{train}} = 100$.  We use a uniform triangular mesh with $2048$ linear elements, and the backward Euler time-stepping scheme with $\Delta t = 1\times 10^{-2}$ for $N_t = 200$ steps until $T = 2$.

\begin{figure}[!ht]
  \centering
  \includegraphics{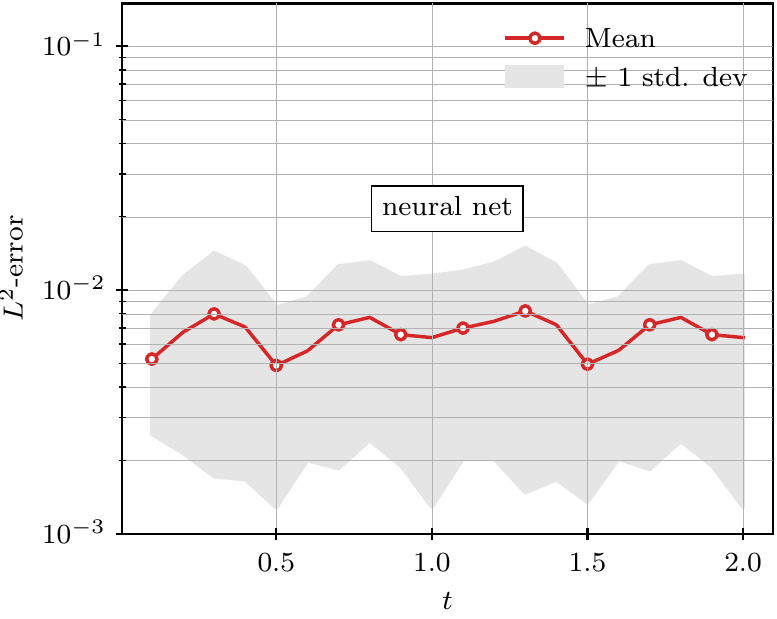}
\caption{\textbf{Non-affine diffusion-reaction}:
Mean relative $L^2$-error with respect to full-order solution. The shaded region denotes $\pm$ one standard deviation over the test set. Timesteps are marked every 10 steps.}\label{fig:ex_3_mean_err}
\end{figure}

As before, we use only the data corresponding to every 10 timesteps, and perform the POD process on $\valpha(t_k;\vmu)$ with $\epsilon_t = \epsilon_{\vmu} = 10^{-5}$, resulting in 7 $L^2$-orthonormal basis functions.  The network is tested on a set of $\Ntrain = 50$ parameters generated with Latin hypercube sampling.  The mean relative error over the test set is shown in~\Cref{fig:ex_3_mean_err}.
We observe that the neural network yields high approximation quality over the test set, with the error for the majority of the samples bounded below $1.6\%$.

\section{Conclusion}

In this work we have developed a method to compute coefficients for
reduced-order models of time-dependent problems in cases of parametric
non-separability.  By training a neural network to identify the evolution of
reduced-basis coefficients, the method is able to accurately capture temporal
behavior in several test cases, including problems that are not parametrically
separable.  Overall, the method provides a flexible approach for developing
robust reduced-order models, while remaining computationally efficient.

There are several aspects of the work that warrant further study, but are
beyond the scope of the immediate paper.  Throughout the work and in related
methods, the neural network architectures and settings have only loosely
been explored.  Here, we have selected parameters in the network, such as the layer width
and the number hidden layers per block, to provide a robust model.  One
benefit of our work is that these are not highly tuned. However, there may be
additional benefits in conducting a deeper parameters study in terms of
computational efficiency. Likewise, we have used modest thresholds in the POD
constructions in the algorithm; a detailed analysis of these parameters and the
associated error estimates would compliment the work presented here.
Finally, it is important to note that the time-based PDEs considered in this
study were selected to highlight the efficacy of the proposed approach.  Yet we
recognize that more complex flows, for example with shocks or developing
boundary layers, may require additional attention.  This would be a natural
direction for continued work.

\section*{Acknowledgments}

The work of R. Patel, supported by the U.S. Department of Energy, Office of Advanced Scientific Computing Research under the Collaboratory on Mathematics and Physics-Informed Learning Machines for Multiscale and Multiphysics Problems (PhILMs) project.  R. Patel and E. C. Cyr are both supported by the Department of Energy early career program. Sandia National Laboratories is a multimission laboratory managed and operated by National Technology and Engineering Solutions of Sandia, LLC, a wholly owned subsidiary of Honeywell International, Inc., for the U.S. Department of Energy’s National Nuclear Security Administration under contract DE-NA0003525. This paper describes objective technical results and analysis. Any subjective views or opinions that might be expressed in the paper do not necessarily represent the views of the U.S. Department of Energy or the United States Government. SAND Number: SAND2021-12950 O

\bibliographystyle{elsarticle-num-names}
\bibliography{refs-nnrb}

\appendix

\section{Notation}\label{sec:notation}

For reference,~\cref{tab:notation} provides a summary of symbols and parameters used throughout the paper.
\begin{table}
  \centering
\begin{tabular}{l l}
  \toprule
  Variable & Definition\\
  \midrule
  $\Omega \in \mathbb{R}^d$ & spatial domain\\
  $T$ & end time of the simulation\\
  $u(\vx, t; \vmu)$ & PDE solution\\
  $F$ & differential operator\\
  $\vmu$ & vector of $P$ parameters, $\mu_1$, \dots, $\mu_P$\\
  $V^h$: $\phi_1$, \dots, $\phi_{N_h}$ & full order basis on mesh $h$\\
  $V^{\rb}$: $\psi_1$, \dots, $\psi_{N_{\rb}}$ & reduced basis on mesh $h$\\
  $t_k$ & $k=1$,\dots,$N_t$ time step $k$ with $t_{N_{t}} = T$\\
  $N_t$ & number of time steps to the end time\\
  $\mathcal{D}$ & Parameter space\\
  $N_\text{s}$ & Number of POD samples\\
  $\mathcal{D}_{\text{POD}}$ & POD sample space\\
  $S_j$ & Data matrix for parameter sample ${\vmu}_j$\\
  $U$ & Data matrix combining all compressed samples\\
  $N_{\rb_j}$ & first stage basis functions for parameter sample $j$ (SVD)\\
  $m_j$ & \# of singular values for parameter sample $j$\\
  $m$ & \# of singular values in the second stage\\
  $\epsilon_{\vmu}$ & stopping criteria, first stage\\
  $\epsilon_{t}$ & stopping criteria, second stage\\
  $\valpha(t;\vmu)$ & vector of full-order basis coefficients\\
  $\vc(t;\vmu)$ & vector of reduced basis coefficients\\
  $u^{rb}$ & reduced order PDE solution\\
  $\mathcal{N}_R$ & mapping of a parameter to RB coefficient\\
  $L_T$ & Least squares single-step loss\\
  $L_T^m$ & Least squares multi-step loss\\
  $\mathcal{T}$ & A training set specifying time nodes and parameters \\
  $\mathcal{B}$ & A mini-batch defined as a subset of $\mathcal{T}$ \\
  $\mathcal{N}_T$ & the neural net and mapping\\
  $\mathcal{D}_{\text{train}}$ & training set of parameters\\
  $N_{\text{train}}$ & number of training samples\\
  $\mathcal{F}$ & Training data features\\
  $\mathcal{T}$ & Training data targets\\
  $L$ & number of layers in a ResNet \\
  $R_i$ & ResNet block\\
  $I$ & Identity\\
  $A^{(j)}$ & Dense NN layer\\
  $W^{(j)}$ & Weight matrix\\
  $b^{(j)}$ & Bias vector\\
  $\sigma$ & activation function\\
  \bottomrule
\end{tabular}
\caption{Notation used throughout.}\label{tab:notation}
\end{table}

\section{Learning Parameters}

In~\cref{table:allexamples} we summarize the
parameters used in each numerical experiment.
\begin{table}
  \centering
  \begin{tabular}{llll}
\toprule
Parameter               & 1D Adv. Diff. (\labelcref{eq:example-1})
                        & 2D Adv. Diff. (\labelcref{eq:example-2})
                        & 2D non-affine (\labelcref{eq:nonaffine})\\
\midrule
\# of residual blocks   & 2                    & 2                    & 1      \\
Hidden layers per block & 2                    & 2                    & 3      \\
Width of each block     & 10, 8                 & 8, 7                 & 10     \\
$m$                     & 2                    & 2                    & 4      \\
Activation function     & elu                  & elu                  & elu    \\
Optimizer               & L-BFGS               & L-BFGS               & L-BFGS \\
Learning rate           & 1.0                  & 1.0                  & 1.0    \\
Max epochs              & 1,250                & 1,250                & 5,000  \\
\multirow{2}{*}{Training set size}
                        & 81 $\vmu$ values     & 50 $\vmu$ values     & 100 $\vmu$ values\\
                        & 100 timesteps         & 50 timesteps         & 20 timesteps    \\
\bottomrule
  \end{tabular}
  \caption{Hyperparameters for the examples in~\cref{sec:numerics}.}\label{table:allexamples}
\end{table}

\end{document}